\documentclass[12pt]{amsart}

\usepackage{amsgen,amsmath,amstext,amsbsy,amsopn,amsfonts,amssymb}
\usepackage{amsthm}
\usepackage[dvips]{graphicx}
\usepackage{epsfig}
\usepackage{ifthen}
\usepackage{amscd}
\usepackage[all]{xy}

\newfont{\fnt}{bbold10 scaled 1250}

\newcommand{\Oo}{\mathcal O}

\newcommand{\Ee}{\mathcal E}
\newcommand{\Ff}{\mathcal F}
\newcommand{\Mm}{\mathfrak M}

\newcommand{\Rr}{\mathbb R}
\newcommand{\Cc}{\mathbb C}

\newcommand{\Zz}{\mathbb Z} 

\newcommand{\Id}{\mbox{\fnt 1}}

\newcommand{\CP}{\mathbb{P}}
 
\newcommand{\Aa}{\mathcal A}
\newtheorem{teor}{Theorem}[section]
\newtheorem{defi}[teor]{Definition}
\newtheorem{coro}[teor]{Corollary}
\newtheorem{lema}[teor]{Lemma}
\newtheorem{prop}[teor]{Proposition}
\newtheorem{conj}[teor]{Conjecture}

\title[Blowups of surfaces and moduli of vector
bundles]{Blowups of surfaces and moduli of holomorphic vector bundles}

\author{Jo\~ao Paulo Santos}

\begin{document}

\bibliographystyle{amsplain}

\setlength{\parindent}{0em}
\setlength{\parskip}{1.5ex}
\newdir{ >}{{}*!/-5pt/\dir{>}}
\newdir{ |}{{}*!/-5pt/\dir{|}}

\begin{abstract}
We examine the moduli of framed 
holomorphic bundles over the blowup of a complex surface, by studying
a filtration induced by the behavior of the bundles on a neighborhood
of the exceptional divisor.
\end{abstract}

\maketitle

\section{Introduction}

The motivation for this paper comes from the study of moduli spaces of
based instantons. In \cite{Buc93}, \cite{Mat00}, it was
shown that the moduli space of based instantons over
a connected sum of $q$ copies of $\CP^2$ is isomorphic as a real analytic
space, to the moduli space of holomorphic bundles over a blow up
of $\CP^2$ at $q$ points, framed at a rational curve $L_\infty\subset\CP^2$.
From the study of this moduli space
we were led to consider the relationship betwen the moduli space
over an algebraic surface $X$
and the moduli space over its blow up $\tilde X$.
Bundles on the blow up of a complex surface have been studied in
\cite{Bru90}, \cite{FrMo88}, \cite{Gas98}, \cite{Gas00}
\cite{BaGa02}. Our approach is inspired by that in \cite{Buc00}.
This paper builds up on results in
\cite{San01}, extending them to the compactification of the moduli space.

\subsection{Results}

Let $X$ be a smooth algebraic surface and let $C\subset X$
be a curve of positive self-intersection. When $X=\CP^2$ we will always
take $C$ as a rational curve. Fix
an ample divisor $H$ and a polynomial $\bar\delta(n)$ with positive
coefficients.
Let $\Mm_k^{r,ss}(X,H,\bar\delta)$ denote the moduli space of 
$(H,\bar\delta)$ semi-stable
pairs $(\Ee,\phi)$ where $\Ee\to X$ is a rank $r$ coherent sheaf with
$c_1(\Ee)=0$, $c_2(\Ee)=k$, and $\phi:\Ee\to\Oo_{C_\infty}^r$
is a non zero homomorphism (the framing). 
Stability of the pair $(\Ee,\phi)$ means the following: 

\begin{defi}\label{defistability}
A pair $(\Ee,\phi)$ is said to be (semi)stable 
with respect to $(H,\bar\delta)$
if for all subsheaves $\Aa\subset\Ee$ we have
\[
{\rm rk\,}\Ee\left(\ 
\chi(\Aa(n))-\varepsilon\bar\delta(n)\ \right)(\leq)<
{\rm rk\,}\Aa\left(\ 
\chi(\Ee(n))-\bar\delta(n)\ \right)
\]
where $\varepsilon=0$ if $\Aa\subset{\rm ker}\,\phi$ and $\varepsilon=1$
otherwise.
\end{defi}

This is a special case of the construction of moduli of
framed sheaves in \cite{HuLe95b}, \cite{HuLe95a}. 
We will omit reference to $(H,\bar\delta)$
unless the dependance on the polarization is important.

In this paper we will be looking at the subspace 
$\Mm_k^r(X)\subset\Mm_k^{r,ss}(X)$ of pairs $(\Ee,\phi)$ where
$\Ee\to X$ is a holomorphic rank $r$ vector bundle,
trivial when restricted to $C$, and 
$\phi$ induces a trivialization. 
Let $\overline{\Mm_k^r(X)}$ denote the closure of $\Mm_k^r(X)$.
The objective of this paper is to prove the theorem

\begin{teor}\label{mainteor}
Let $\pi:\tilde X\to X$ be the blow up of $X$ at a point $x_0\notin C$.
Given a sheaf $\Ee\to\tilde X$, let $\pi_*^{\vee\vee}\Ee\to X$ 
be the sheaf defined
by $\pi_*^{\vee\vee}\Ee(U)=\pi_*\Ee(U\setminus\{x_0\})$.
\begin{enumerate}
\item Let
\[
S_i\Mm_k^r(\tilde X)=\left\{\,(\Ee,\phi)\,|\,
c_2\left(\pi_*^{\vee\vee}\Ee\right)=i\,\right\}
\]
Then the map $\pi_*^{\vee\vee}:S_i\Mm_k^r(\tilde X)\to\Mm_i^r(X)$
is a topologically trivial fibration with fiber $S_0\Mm_{k-i}(\tilde\CP^2)$;
\item Let 
\[
F_i\overline{\Mm_k(\tilde X)}=\left\{\,(\Ee,\phi)\,|\,
c_2\left(\pi_*^{\vee\vee}\Ee\right)\leq i\,\right\}
\]
and let 
$S_i\overline{\Mm_k}=F_i\overline{\Mm_k}\setminus F_{i-1}\overline{\Mm_k}$.
Then the map $\pi_*^{\vee\vee}$ extends to a map
$\pi_\bullet:F_i\overline{\Mm_k(\tilde X)}\to\overline{\Mm_i(X)}$.
The restriction of this map to $S_i\overline{\Mm_k(\tilde X)}$
is a fibration with fiber $S_0\overline{\Mm_{k-i}(\tilde\CP^2)}$.
\end{enumerate}
\end{teor}

The filtration $F_i\overline{\Mm_k^r}$ induces a
spectral sequence converging to 
$H_*\left(\overline{\Mm_k^r}\right)$ with $E^1$ term given by
\[
E^1_{p,q}=H_{p+q}\left(\,
F_p\overline{\Mm_k^r(\tilde X)},\,
F_{p-1}\overline{\Mm_k^r(\tilde X)}\,\right)
\]
As a corolary we will prove

\begin{coro}\label{maincoro}
Let $\hat\Mm_i(X)=\pi_\bullet F_{i-1}\overline{\Mm_k^r(\tilde X)}\subset
\overline{\Mm_i(X)}$. Then
there is a spectral sequence converging to
$H_*\left(F_i\overline{\Mm},F_{i-1}\overline{\Mm}\right)$
with $E^2$ term given by
\[
E^2_{p,q}=H_p\left(\,\overline{\Mm_i(X)},\,
\hat\Mm_i(X);\,\,H_q(S_0\overline{\Mm_{k-i}(\tilde\CP^2)})\,\right)
\]
\end{coro}

We conjecture a similar result holds in the non compactified case..
Here the triviality of the bundle would lead to a more powerful result:

\begin{conj}
There is an isomorphism
\[
H_*\left(F_i\Mm_k^r(\tilde X),F_{i-1}\Mm_k^r(\tilde X)\right)\cong
H_*\left(\pi_\bullet F_i\Mm_k^r(\tilde X),\,\Mm_i^r(X)\,\right)\otimes
H_*\left(S_0\Mm_{k-i}^r(\tilde\CP^2)\right)
\]
\end{conj}

The plan of the paper is as follows:
In section 2 we recall the definitions of the moduli space and prove
some theorems concerning stability;
In section 3 we prove the second part of the theorem \ref{mainteor},
and corollary \ref{maincoro}.
In section 4 we prove the first part of the theorem \ref{mainteor}.
Finally in section 5 we give a monad description of the space 
$S_0\Mm_k(\tilde\CP^2)$.

\section{Definition of the moduli space}

The objective of this section is to recall the definitions of the moduli
spaces we will be considering, and to prove some stability results.

First we introduce a technical restriction on the pairs
$(X,C)$ we will be considering in this paper.
Then we present two definitions of the moduli space: 
Following \cite{Lub93}, we
look at $\Mm_k^r$ as the space of holomorphic structures on a fixed
topological bundle $E_{top}$. The second definition uses the language
of moduli functors, following \cite{HuLe95b}, \cite{HuLe95a}.

\subsection{Admissible pairs}

\begin{defi}
Let $X$ be an algebraic surface and let $C\subset X$ be a divisor
with $C^2>0$. We say the pair $(X,C)$ is admissible if there is
an ample divisor $H$ such that the set 
$\{\,c_1(\Aa)\cdot H\,\}_{\Aa\in\mathcal C_k^r}$ is bounded above
for all $k,r$.
Here $\mathcal C_k^r$ denotes the set of sheaves $\Aa\to X$ such that
\begin{itemize}
\item $c_1(\Aa)\cdot C\leq0$;
\item There is a torsion free rank $r$
sheaf $\Ee\to X$ with $c_1(\Ee)=0$, $c_2(\Ee)=k$ such that $\Aa\subset\Ee$.
\end{itemize}
\end{defi}

Clearly, if $C$ is ample, $(X,C)$ is admissible.

\begin{prop}
Let $(X,C)$ be an admissible pair with ample divisor $H$
and let $\tilde X\to X$ be the blow
up of $X$ at a point $x_0\notin C$. Then the pair
$(\tilde X,C)$ is admissible with respect to the ample divisor
$\tilde H=H+C-L$ where $L$ denotes the exceptional divisor..
\end{prop}
\begin{proof}
Let $\Aa\in\mathcal C_k^r$, $\Aa\subset\Ee$.
Then $(\pi_*\Aa)^{\vee\vee}\subset(\pi_*\Ee)^{\vee\vee}$. Hence,
the proposition will follow if we prove that
$-c_1(\Aa)\cdot L\leq 2k$.

We may assume without loss of generality that $\Aa,\Ee$ are locally free.
Let $\mathcal T={\rm Tor\,}\Ee/\Aa$. Then we have the diagram
\[
\xymatrix{
&&&0&0\\
&&&\mathcal Q\ar[u]\ar[ur]\\
0\ar[r]&\Aa\ar[r]&\Ee\ar[r]\ar[ur]&\Ee/\Aa\ar[r]\ar[u]&0\\
&\mathcal K\ar[ur]&&\mathcal T\ar[u]\\
0\ar[ur]&&&0\ar[u]
}\]
for some sheaves $\mathcal Q$ and $\mathcal K$. Now observe that
\[
c_1(\mathcal K)\cdot H\geq c_1(\Aa)\cdot H
\]
So we may assume without loss of generality that the quotient
$\Ee/\Aa$ is torsion free. Now write
\[
\Aa|_{L}=\sum_j\Oo(a_j)\,,\ 
\Ee|_{L}=\sum_l\Oo(b_l)
\]
Then $a_j=b_{l_j}$ so
\[
|c_1(\mathcal A)\cdot L|\leq \sum_j|a_j|\leq\sum_l|b_l|\leq 2c_2(\Ee)
\]
which concludes the proof.
\end{proof}

\subsection{Analytic definition}

Let $E\to X$ be an $SU(r)$ topological bundle with $c_2(E)=k$.

Let $\mathcal C(X,E)$ be the space of pairs $(\bar\partial,\phi)$ where
$\bar\partial:\Omega^0(E)\to\Omega^{0,1}(E)$ 
is a holomorphic structure on $E$ holomorphically
trivial on $C$ and
$\phi:E|_{C}\to\Oo_C^r$ is an isomorphism of holomorphic bundles.

\begin{prop}
The group $Aut(E)$ acts freely on ${\mathcal C}(X,E)$ and the quotient has the 
structure of a finite dimensional Hausdorff complex analytic space.
\end{prop}
\begin{proof}
It follows from
\[
\forall_{k>0}H^0({\rm End\,}\Ee_0\otimes\Oo_{C}(-k))=0
\]
where $\Ee_0=\Oo_C^r$ (see \cite{Lub93} theorem 1.1 and lemma 2.6).
\end{proof}

\begin{defi}
We define $\Mm(X,E)$
as the quotient $\Cc(X,E)/Aut(E)$.
We will also use the notation
$\Mm_k^r(X)$.
\end{defi}

\subsection{Algebraic definition}

Now we present the algebraic definition. For details see
\cite{HuLe95b}, \cite{HuLe95a}.

Let $\Ee_0=\Oo_C^r$. A family of framed sheaves
parametrized by a Noetherian scheme $T$ consists of a pair $(\Ff,\alpha)$
where $\Ff$ is a coherent $\Oo_{T\times X}$ module, flat over $T$, and
$\alpha:\Ff\to \Oo_T\otimes\Ee_0$ is a homomorphism with $\alpha_t\neq0$.
A homomorphism of families $(\Ff,\alpha)\to(\Ff',\alpha')$
is a homomorphism of sheaves
$\Ff\to\Ff'$ compatible with $\alpha,\alpha'$.

The moduli functor 
$\underline{\mathcal M^{ss}(X,\Ee_0)}$ is the functor
from $(Schemes)$ to $(Sets)$ 
that to a scheme $T$ associates the set
of isomorphism classes of flat families of semistable pairs parametrized by $T$
(recall the definition of stability \ref{defistability}).
In a similar way we define the functor $\underline{\mathcal M^{s}(X,\Ee_0)}$
by replacing the word semistable with stable in the definition.

In \cite{HuLe95a} Huybrechts and Lehn proved

\begin{teor}
There is a projective scheme $\Mm_k^{r,ss}(X)$ which corepresents
the functor $\underline{\mathcal M^{ss}(X)}$. 
Moreover there is an open subscheme
$\Mm_k^{r,s}$ which is a fine moduli space representing
$\underline{\mathcal M^s(X)}$.
\end{teor}

$\Mm_k^r(X)\subset\Mm_k^{r,ss}(X)$ 
is defined as the subspace of pairs $(\Ee,\phi)$ where 
$\Ee\to X$ is a holomorphic vector bundle
trivial when restricted to $C$, and 
$\phi$ induces a trivialization. 
$\overline{\Mm_k^r(X)}$ is defined as the closure of $\Mm_k^r(X)$.

\subsection{Some results about stability}

A pair $(\Ee,\phi)$ is said to be $\mu$-stable with respect to 
$(H,\delta)$ if for every subsheaf $\Aa\subset\Ee$ we have
${\rm rk\,\Ee}(c_1(\Aa)\cdot H-\varepsilon\delta)\leq
-{\rm rk\,\Aa}\,\delta$ where $\varepsilon$ is defined as 
in \ref{defistability}.
Then, from Riemann-Roch we get the implications
\[
\mu-\mbox{stable}\Rightarrow\mbox{stable}\Rightarrow\mbox{semistable}
\Rightarrow\mu-\mbox{semistable}
\]

\begin{prop}
Let $H$ be an ample divisor and consider a polarization
$(H+MC,\delta'+n\delta)$ with $M>\delta$. Then, if a pair
$(\Ee,\phi)$ is semistable, $\Ee$ is torsion free.
\end{prop}
\begin{proof}
We apply the semistability 
condition to $\Aa={\rm Tor}\,\Ee$. We want to show $\Aa=0$.
We divide the proof into two steps:
\begin{enumerate}
\item $\mu$-semistability implies
$c_1(\Aa)\cdot(H+MC)\leq\varepsilon\delta$.
Suppose $c_1(\Aa)\neq0$. Then $\varepsilon=1$
so the restriction of $\phi$ to $\Aa$ is not identically zero.
But then we must have $c_1(\Aa)\cdot C>0$ so
$(H+MC)\cdot c_1(\Aa)\geq M>\delta$ contradicting semistability.
Hence $c_1(\Aa)=0$.
\item Since $c_1(\Aa)=0$,
$\Aa$ is supported in codimension 2, so $\varepsilon=0$ and $c_2(\Aa)\leq0$.
On the other hand semistability implies
$\chi(\Aa)=-c_2(\Aa)\geq0$. Hence
$c_2(\Aa)=0$. So $\Aa=0$.
\end{enumerate}
\end{proof}

\begin{lema}\label{lemastability}
Let $(X,C)$ be an admissible pair with ample divisor $H$.
Fix a pair $(\Ee,\phi)$. Then there is an integer
$\delta_0>0$ such that, for all $\delta>\delta_0$ there is an integer $M_0$
depending on $\delta$ such that for all $M>M_0$ the following holds
with respect to the choice of polarization 
$(H+MC,\delta'+\delta n)$:
\begin{enumerate}
\item Let $\Aa\subset\Ee$ be such that  $c_1(\Aa)\cdot C<0$.
Then $\Aa$ is not destabilizing.
\item Let $\Aa\subset\Ee$ be such that $c_1(\Aa)\cdot C\leq0$ and
$\varepsilon=1$. Then $\Aa$ is not destabilizing.
\end{enumerate}
\end{lema}
\begin{proof}
For simplicity we assume $c_1(\Ee)=0$. 
Let $r_\Ee={\rm rk}\,\Ee,\,r_\Aa={\rm rk}\,\Aa$.
Since $(X,C)$ is admissible
we can choose $\delta_0>r_\Ee c_1(\Aa)\cdot H$.
Now pick any $\delta>\delta_0$. Then choose $M_0$ so that
$M_0>c_1(\Aa)\cdot H+\delta$. Finally pick any $M>M_0$. Let $H_M=H+MC$.
\begin{enumerate}
\item Assume $c_1(\Aa)\cdot C<0$. Then
\[
r_\Ee c_1(\Aa)\cdot H_M\leq r_\Ee(c_1(\Aa)\cdot H-M)<
-r_\Ee\delta<-r_\Aa\delta.
\]
So, we get
\[
\frac{c_1(\Aa)\cdot H_M-\varepsilon\delta}{r_\Aa}<
\frac{c_1(\Ee)\cdot H_M-\delta}{r_\Ee}
\]
\item Assume $c_1(\Aa)\cdot C\leq0$. Then
\[
r_\Ee(c_1(\Aa)\cdot H_M-\delta)\leq r_\Ee c_1(\Aa)\cdot H-
r_\Ee\delta<(1-r_\Ee)\delta
\leq -r_\Aa\delta.
\]
Now, since by assumption $\varepsilon=1$,
\[
\frac{c_1(\Aa)\cdot H_M-\varepsilon\delta}{r_\Aa}<
\frac{c_1(\Ee)\cdot H_M-\delta}{r_\Ee}
\]
\end{enumerate}
\end{proof}

As a corollary we have

\begin{prop}\label{propstability}
There is an ample divisor $H$ and a polynomial $\bar\delta$ such that,
for any $M>0$, the following holds
with respect to the polarization $(H+MC,\bar\delta)$:

Let $\Ee\to X$ be a torsion free sheaf,
and let $\phi:\Ee\to\Ee_0$ induce an isomorphism $\Ee|_C\to\Oo_C^r$.
Then the pair $(\Ee,\phi)$ is stable.
\end{prop}
\begin{proof}
Let $(H+MC,\delta'+n\delta)$ be such that the conclusions of 
lemma \ref{lemastability} hold.
The proof will follow from the following
statements:
\begin{itemize}
\item If $\Aa\subset\Ee$ then $C\cdot c_1(\Aa)\leq 0$.
\item If $\varepsilon=0$ then $C\cdot c_1(\Aa)<0$.
\end{itemize}
The first statement follows from 
$\Aa|_{C}\subset\Ee|_{C}\cong\Oo^r$. To prove the second
we observe that ${\rm Ker}\,\phi=\Ee(-C)$ so if 
$\Aa\subset{\rm Ker}\,\phi\subset\Ee$ 
then $\Aa(C)\subset\Ee$. The result follows.
\end{proof}

Now we prove a converse to lemma \ref{lemastability}:

\begin{prop}\label{propstability2}
Let $(H,\bar\delta(n)=\delta'+n\delta)$ 
be a polarization for which the conclusions of lemma
\ref{lemastability} hold.
Let $(\Ee,\phi)\in\overline{\Mm_k^r(X)}$. Then, for all subsheaves
$\Aa\subset\Ee$, we have either $c_1(\Aa)\cdot C<0$ or
$c_1(\Aa)\cdot C=0,\varepsilon=1$.
\end{prop}
\begin{proof}
Let $(\Ee,\phi)\in\overline{\Mm_k^r(X)}\subset\Mm^{ss}(H,\bar\delta)$ 
and suppose there was some
$\Aa\subset\Ee$ with either $c_1(\Aa)\cdot C>0$ or
$c_1(\Aa)\cdot C=\varepsilon=0$. We want to show this is not possible.
Consider a new polarization
$(H_M,\bar\delta_M)$ with $H_M=H+M_HC$, 
$\bar\delta_M(n)=\bar\delta(n)+nM_\delta$, for some constants 
$M_H>M_\delta>0$.
Then
\begin{enumerate}
\item We claim that
if $(\Ee',\phi')$ is any pair semistable with respect to
$(H_M,\bar\delta_M)$, then $(\Ee',\phi')$ is stable with respect to 
$(H,\bar\delta)$.
Let $\Aa'\subset\Ee'$. By lemma \ref{lemastability}
we may assume either $c_1(\Aa')\cdot C>0$ or $c_1(\Aa')\cdot C=\varepsilon=0$.
Then $\mu$-semistability with respect to $(H_M,\bar\delta_M)$ implies
\begin{multline*}
0\geq r_\Ee(c_1(\Aa')\cdot H_M-\varepsilon(\delta+M_\delta))+
r_\Aa(\delta+M_\delta)\geq\\
\geq r_\Ee(c_1(\Aa')\cdot H-\varepsilon\delta)+r_\Aa\delta
+M_Hr_\Ee c_1(\Aa')\cdot C+M_\delta(1-r_\Ee)>\\
>r_\Ee(c_1(\Aa')\cdot H-\varepsilon\delta)+r_\Aa\delta
\end{multline*}
where $r_\Ee={\rm rk\,}\Ee,\,r_\Aa={\rm rk\,}\Aa$. So
$(\Ee',\phi,)$ is $\mu$-stable with respect to $(H,\bar\delta)$.
\item We claim that for $M_H>M_\delta\gg0$,
$(\Ee,\phi)$ is unstable with respect to $(H_M,\bar\delta_M)$;
We have two cases: 
\begin{itemize}
\item If $c_1(\Aa)\cdot C>0$ then 
we just have to choose $M_H$ big enough so that $(\Ee,\phi)$ is
$\mu-$unstable;
\item If $c_1(\Aa)\cdot C=\varepsilon=0$ then stability questions are not
affected by $M_H$ and making $M_\delta$ big enough we can make
$(\Ee,\phi)$ $\mu-$unstable.
\end{itemize}
\end{enumerate}
Using the fact $\Mm^{ss}$ corepresents 
the moduli functor, property (1) implies there is a map 
$\Mm^{ss}(H_M,\bar\delta_M)\to\Mm^{s}(H,\delta)$. Let
$D\subset\overline{\Mm_k^r(X)}\subset\Mm^{ss}(H,\bar\delta)$ be a one 
dimensional disk and suppose the origin $0\in D$ corresponds to
$(\Ee,\phi)$ and $D\setminus0$ is contained in $\Mm_k^r(X)$.
Then, by proposition \ref{propstability}, 
the restriction of the universal family to $D\setminus 0$ gives a family
of $(H_M,\bar\delta_M)$-stable pairs. Hence we have a map
$D\setminus0\to\Mm^{ss}(H_M,\bar\delta_M)$. Since
$\Mm^{ss}(H_M,\bar\delta_M)$ is projective this map extends to a map
$D\to\Mm^{ss}(H_M,\bar\delta_M)$. We get a commutative diagram
\[\xymatrix{
D\ar[r]\ar[dr]&\Mm^{ss}(H_M,\bar\delta_M)\ar[d]\\
&\Mm^{ss}(H,\bar\delta)
}\]
But this contradicts the fact that 
$(\Ee,\phi)\notin\Mm^{ss}(H_M,\bar\delta_M)$. This concludes the proof.
\end{proof}

As a corollary we have

\begin{coro}
For a convenient choice of polarization,
$\overline{\Mm_k^r(X)}\subset\Mm_k^{r,s}$.
\end{coro}

\section{The map $\pi_\bullet$}

In this section we prove the second part of theorem \ref{mainteor}.
Recall that
\begin{align*}
&F_i\overline{\Mm_k(\tilde X)}=\left\{\,(\Ee,\phi)\,|\,
c_2\left(\pi_*^{\vee\vee}\Ee\right)\leq i\,\right\}\\
&S_i\overline{\Mm_k(\tilde X)}=\left\{\,(\Ee,\phi)\,|\,
c_2\left(\pi_*^{\vee\vee}\Ee\right)= i\,\right\}
\end{align*}
First we want to show the existence of a map
$\pi_\bullet:F_i\overline{\Mm_k(\tilde X)}\to\overline{\Mm_i(X)}$ extending
$\pi_*^{\vee\vee}$. We will define this map locally.
Let $T\subset F_i\overline{\Mm_k^r(\tilde X)}$ be an open subset
and consider the universal family $(\Ff,\alpha)$,
$\Ff\to T\times\tilde X$, $\alpha:\Ff\to p^*\Ee_0$. Let $L$
be the exceptional divisor.

\begin{prop}\label{prop31}
For each $N\in\Zz$ define the family $\pi_\bullet^N\Ff\to T\times X$ by
\[
\pi_\bullet^N\Ff(W)=(\Id\times \pi)_*\Ff(-NL)(W\setminus(T\times x_0))
\]
For each  $t=[\Ee,\phi]\in T$ let $\imath_t:\{t\}\times X\to T\times X$
be the inclusion.
Then we have, for $N\gg0$,
\begin{enumerate}
\item $(\pi_\bullet^N\Ff,\alpha)$ is a flat family of stable pairs over $X$,
hence it induces a map 
$\pi_\bullet:T\to\overline{\Mm_k^r(X)}$ given by
$t\mapsto (\imath_t^*\pi_\bullet\Ff,\alpha|_t)$.
\item If $t\in T\cap S_i\Mm_k^r$, then 
$\imath_t^*\pi_\bullet\Ff=(\pi_*\Ee)^{\vee\vee}$. Hence $\pi_\bullet$ extends
$\pi_*^{\vee\vee}$.
\end{enumerate}
\end{prop}
\begin{proof}
We begin by showing part (1). We have to show flatness and stability.
\begin{itemize}
\item We prove flatness in two steps. First we show that
$\Ff'=(\Id\times\pi)_*\Ff(-NL)$ is flat. It is enough to show that
the Hilbert polynomial of $\imath^*_t\Ff'$
is constant with $t$ (see \cite{HuLe97}, proposition 2.1.2).
This follows from the Grothendieck-Riemann-Roch theorem applied
to $\imath^*_t(\Id\times\pi)_*\Ff(-NL)=\pi_*\imath_t^*\Ff(-NL)$ 
(for this equality see the appendix) plus
the vanishing of the higher direct image sheaves $R^i\pi_*$ for $N\gg0$.

Now we want to show $\pi^N_\bullet\Ff$ is flat. It is enough to check for
points $(t_0,x_0)$ for some $t_0\in T$.
Fix an ideal $I\subset\Oo_{t_0,T}$.
Let $\sum_ia_i\otimes m_i\in I\otimes\pi_\bullet\Ff_{t_0,x_0}$ and 
assume $\sum_ia_im_i=0\in\pi_\bullet\Ff_{t_0,x_0}$.
We want to show $\sum_ia_i\otimes m_i=0$.

Pick an open set $W\subset T\times X$ such that the following holds:
\begin{itemize}
\item For all $i$, $a_i\in \Oo_T(p(W))$ and
$m_i\in\Ff'(W\setminus (T\times x_0))$;
\item There are generators $f_1,\ldots,f_j$ of $I$ such that
$f_j\in\Oo_T(p(W))$.
\end{itemize}
Then $\sum_ia_im_i=0\in\Ff'(W\setminus (T\times x_0))$ hence
$\sum_ia_im_i=0\in\Ff'_{(t,x)}$ for any $x\neq x_0$ and any $t$.

Define the sheaf $\mathcal I$ as the sheaf of ideals
$\mathcal I(V)=\left\langle\,I\cap\Oo_T(p(W))\,\right\rangle\subset
\Oo_T(V)$. Then $\mathcal I_{t_0}=I$. Then
by flatness of $\Ff'$ it follows that 
$\sum_ia_i\otimes m_i=0\in\mathcal I_t\otimes\Ff'_{(t,x)}$.
Hence 
$\sum_ia_i\otimes m_i=0\in(\mathcal I\otimes\Ff')(W\setminus(T\times x_0))$.
Now
\[
(\mathcal I\otimes\Ff')(W\setminus(T\times x_0))=
\mathcal I(p(W))\otimes\pi_\bullet\Ff(W)
\]
The flatness of $\pi_\bullet\Ff$ follows.
\item Now we prove stability. We want to show that, for every $t\in T$,
$\imath_t^*\pi_\bullet\Ff$ is stable. Let $\Ee\to X$ be a sheaf defined by
$\Ee(V)=\imath_t^*\pi_\bullet\Ff(V\setminus\{x_0\})$. Then stability
of $\Ee$ is equivalent to stability of $\imath_t^*\pi_\bullet\Ff$.
Let $\Aa\subset\Ee$. We may assume $\Aa$ is locally free at $x_0$.
We claim that, for $M\gg0$,  $\pi^*\Aa(-ML)\subset\imath_t^*\Ff$.
To see this notice that $\Ee=\pi_*\imath_t^*\Ff(ML)$ hence we have the 
inclusion map $\pi^*\Aa\to\pi^*\pi_*\imath_t^*\Ff(ML)\to\imath_t^*\Ff(ML)$.

Now notice that $\Aa$ and $\pi^*\Aa(-ML)$ are isomorphic on 
$X\setminus x_0$. Hence,
since $\imath_t^*\Ff$ is stable, it follows by
proposition \ref{propstability2} that either $c_1(\Aa)\cdot C<0$ or
$c_1(\Aa)\cdot C=0,\varepsilon=1$. But then, by lemma \ref{lemastability},
$\Aa$ is not destabilizing.
\end{itemize}
Now we prove statement (2): $\pi_\bullet$  restricted to 
$S_i\Mm_k^r(\tilde X)$ 
is given by $(\Ee,\phi)\mapsto((\pi_*\Ee)^{\vee\vee},\phi)$.
That is, if $t=[\Ee,\phi]\in S_i\Mm_k^r(\tilde X)$ then
$\imath_t^*\pi_\bullet\Ff=(\pi_*\Ee)^{\vee\vee}$.

First observe that, for $t=[\Ee,\phi]\in S_i\Mm_k^r\cap T$, 
$\imath_t^*\pi_\bullet\Ff=\pi_*\Ee$ on $X\setminus\{x_0\}$ hence
$(\imath_t^*\pi_\bullet\Ff)^{\vee\vee}=(\pi_*\Ee)^{\vee\vee}$.
In particular $c_2((\imath_t^*\pi_\bullet\Ff)^{\vee\vee})=i$.
Now flatness implies $c_2(\imath_t^*\pi_\bullet\Ff)$ is constant with $t$
so it is enough to show that for some $t\in S_i\Mm_k^r\cap T$,
$\imath_t^*\pi_\bullet\Ff$ is locally free.

To prove this last statement observe that 
in $(S_i\Mm_k^r\cap T)\times X$ we have
$\pi_\bullet\Ff=((\Id\times\pi)_*\Ff)^{\vee\vee}$. 
Hence its singularities lie in codimension three.  This implies
the desired result.
\end{proof}

Now we want to show that the restriction of $\pi_\bullet$ to
$S_i\overline{\Mm_k(\tilde X)}$ is a fibration with fiber
$S_0\overline{\Mm_{k-i}(\tilde\CP^2)}$.

Let $\check\Mm_i^r(X)$ be the subspace of pairs $(\Ee,\phi)$ such that
$\Ee$ is locally free at $x_0$.

\begin{prop}
$\pi_\bullet^{-1}\check\Mm_i^r(X)=S_i\overline{\Mm_k^r(\tilde X)}$
\end{prop}
\begin{proof}
We begin by
remarking that, for any sheaf $\Ee\to\tilde X$, 
$\pi_\bullet\Ee$ and $\pi_*^{\vee\vee}\Ee$ coincide over $X\setminus\{x_0\}$.
So $(\Ee,\phi)\in\pi_\bullet^{-1}\check\Mm_i^r(X)$ if and only if
$\pi_\bullet\Ee=\pi_*^{\vee\vee}\Ee$. Now
\begin{itemize}
\item Suppose $\pi_\bullet\Ee=\pi_*^{\vee\vee}\Ee$. Then, by definition
of $S_i$, $(\Ee,\phi)\in S_i\overline{\Mm_k^r(\tilde X)}$.
\item If $(\Ee,\phi)\in S_i\overline{\Mm_k^r(\tilde X)}$ then
$c_2(\pi_*^{\vee\vee}\Ee)=c_2(\pi_\bullet\Ee)$ hence 
$\pi_\bullet\Ee=\pi_*^{\vee\vee}\Ee$.
\end{itemize}
\end{proof}

\begin{prop}\label{prop33}
The restriction of $\pi_\bullet$ to $\pi_\bullet^{-1}(\check\Mm_i^r(X))$ 
is a fibration with fiber
$S_0\overline{\Mm_{k-i}^r(\tilde\CP^2)}$.
\end{prop}
\begin{proof}
Let $T_X\subset \check\Mm_i^r(X)$ be an open set. We want to
build an isomorphism 
$T_X\times S_0\overline{\Mm(\tilde\CP^2)}\cong
\pi_\bullet^{-1}(T_X)$. Let 
$T_P\subset\Mm_{k-i}^r(\tilde\CP^2)$ be an open set.
Consider the universal families $(\Ff_X,\alpha_X)$ over $T_X$ and 
$(\Ff_P,\alpha_P)$ over $T_P$. The next step is to build trivializations
$\psi_X,\psi_P$ of $\Ff_X,\Ff_P$:
\begin{itemize}
\item For $T_X$ small enough we can choose
a neighborhood of $x_0$, $U\subset X$, such that $\Ff_X$ is free on
$T_X\times U$. Fix a trivialization
$\psi:\Ff|_{T_X\times U}\to\Oo_{T_X\times U}^r$;
\item By definition
of $S_0\overline{\Mm_{k-i}^r(\tilde\CP^2)}$, for any $t\in T_P$ the sheaves
$\imath_t^*\Ff_P|_{\tilde\CP^2\setminus L}$ are free. It follows that,
for $T_P$ small enough,
$\Ff_P|_{\tilde\CP^2\setminus L}$ is a free sheaf. Then Hartog's theorem
implies there is a unique isomorphism 
$\psi_P:\Ff_P|_{\tilde\CP^2\setminus L}\to\Oo^r_{\tilde\CP^2\setminus L}$
which is compatible with $\alpha_P$.
\end{itemize}
Now consider the sheaf $\Ff_{XP}\to\tilde X\times T_X\times T_P$ given by
\[
\Ff_{XP}=p_X^*\,\Ff_X|_{X\setminus\{x_0\}}\bigcup_{\psi_P^{-1}\psi_X}
p_P^*\,\Ff_P|_{\tilde\CP^2\setminus E}
\]
where  $p_X:T_X\times T_P\to T_X$, $p_P:T_X\times T_P\to T_P$ are the 
projections, and let $\alpha_{XP}=\alpha_Xp_X$. We claim 
$(\Ff_{XP},\alpha_{XP})$ is a flat family of stable framed sheaves. 
Flatness follows since both $p_X^*\,\Ff_X$ and $p_P^*\,\Ff_P$
are flat (see \cite{Har77}, proposition III.9.2). To prove stability
let $(t_X,t_P)\in T_X\times T_P$ and let $\Aa\subset\imath_t^*\Ff_{XP}$.
Then $(\pi_*\Aa)^{\vee\vee}\subset\imath_{t_X}^*\Ff_X$.
Now we repeat the argument used in the proof of proposition \ref{prop31}.

The family $(\Ff_{XP},\alpha_{XP})$ induces a map
$g_{\psi}:T_X\times T_P\to\Mm_k^r(\tilde X)$ such that 
$\pi_\bullet g_\psi=p_X$. 
We want to build an inverse
to this map, $\pi_\bullet\times g_{\tilde\psi}$, where
\[
g_{\tilde\psi}:\pi_\bullet^{-1}(T_X)\to
S_0\overline{\Mm_{k-i}^r(\tilde\CP^2)}
\]
Let $T\subset\pi_\bullet^{-1}(T_X)$ and consider the universal family
$(\Ff,\alpha)$ over $T$. We have an isomorphism
\[
\tilde\psi:\Ff|_{T\times\tilde U\setminus L}=
\pi_\bullet\Ff|_{T\times U\setminus\{x_0\}}
\cong(\Id\times\pi_\bullet)^*\Ff_X|_{T\times U\setminus\{x_0\}}
\overset{(\Id\times\pi_\bullet)^*\psi}{\longrightarrow}
\Oo^r_{T\times U\setminus\{x_0\}}
\]
Then we define the sheaf over $\tilde\CP^2\times T$
\[
\Ff_p=\Ff|_{\tilde U}\bigcup_{\tilde\psi}
\Oo^r_{\CP^2\setminus\{x_0\}}
\]
As above, this is a flat family of framed sheaves inducing the 
desired map $g_{\tilde\psi}$.

Now it is a direct verification
to check that 
$\pi_\bullet\times g_{\tilde\psi}$ is the inverse map of $g_\psi$.
\end{proof}

Now we turn to the proof of corollary \ref{maincoro}. We will need the lemma:

\begin{lema}\label{lemaVW}
Let $f:X\to Y$ be a proper map betwen metric spaces $X,Y$.
Let $C\subset Y$ be a closed subspace.
Then, for any neighborhood $W$ of $f^{-1}(C)$,
there is a neighborhood $V$ of $C$ such that 
$f^{-1}(V)\subset W$.
\end{lema}
\begin{proof}
We claim that,
for any $y\in C$ we can build a neighborhood $V_y$ of $y$ such that
$f^{-1}(V_y)\subset W$: if not we could build a sequence
$x_n\notin W$ with $f(x_n)\to y$. Then properness of $f$ leads 
to a contradiction. Now, just take $V=\bigcup_{y\in C}V_y$.
\end{proof}

\begin{proof}[Proof of Corollary \ref{maincoro}]
To simplify notation we will write
$F_i=F_i\overline{\Mm(\tilde X)}$, $\hat\Mm_i=\hat\Mm_i(X)$.
Recall that $\hat\Mm_i=\pi_\bullet F_{i-1}$.
We begin by building open sets $V_0,V_1\subset\overline{\Mm_i^r(X)}$
and $W_0,W_1\subset F_i\overline{\Mm_k^r(\tilde X)}$
such that 
\begin{enumerate}
\item $\pi_\bullet^{-1}V_1\subset W_1\subset 
\pi_\bullet^{-1}V_0\subset W_0$;
\item $\hat\Mm_i\subset V_j, j=0,1$ are strong deformation retracts and
$V_1\setminus\hat\Mm_i\to V_0\setminus\hat\Mm_i$ is a homotopy equivalence;
\item $F_{i-1}\subset W_n,n=0,1$ are strong deformation
retracts and $W_1\setminus F_{i-1}\to W_0\setminus F_{i-1}$ is a homotopy
equivalence.
\end{enumerate}
The existence of such neighborhoods follows from lemma \ref{lemaVW}
and \cite{Hir75}.

Now it follows from the five lemma applied
to the homotopy exact sequence 
coming from the fibrations $\pi^{-1}_\bullet V_j\setminus F_{i-1}\to
V_j\setminus\hat\Mm_i$, $j=1,2$, that the inclusion
$\pi_\bullet^{-1}V_1\setminus F_{i-1}\to
\pi_\bullet^{-1}V_0\setminus F_{i-1}$ is a homotopy equivalence.
Now, using the inclusionmaps in (1) we see that
the spaces $\pi_\bullet^{-1}V_1\setminus F_{i-1}$ and
$W_1\setminus F_{i-1}$ are homotopycally equivalent.
Now, by excision,
\[
H_*(F_i,F_{i-1})\cong H_*(S_i,W_1\setminus F_{i-1})\cong
H_*(S_i,\pi_\bullet^{-1}V_1\setminus F_{i-1})
\]
where $S_i=S_i\overline{\Mm_k^r(\tilde X)}$.
To conclude the proof we apply the relative Leray-Serre spectral sequence to
the pair of fibrations $S_i\to\overline{\Mm_i}\setminus\hat\Mm_i$ and 
$\pi_\bullet^{-1}V_1\setminus F_{i-1})\to V_1\setminus\hat\Mm_i$.
The $E^2$ term is
\[
E^2_{p,q}=H_p\left(\,\overline{\Mm_i}\setminus\hat\Mm_i,V_1\setminus\hat\Mm_i;
H_q(S_0\overline{\Mm_{k-i}(\tilde\CP^2)})\,\right)
\]
To finish the proof we apply excision.
\end{proof}

\section{Stratification of $\Mm_k^r(\tilde X)$}

The results in this section first appeared in \cite{San01}.
Recall that
\[
S_i\Mm_k^r(\tilde X)=\left\{(\Ee,\phi)\in\Mm_k^r(\tilde X)\,|\,
c_2\left(\,(\pi_*\Ee)^{\vee\vee}\,\right)=i\right\}
\]
The objective of this section is to prove part one of theorem
\ref{mainteor}:

\begin{teor}\label{teor51}
The map
\[
\pi_\bullet:S_i\Mm_k^r(\tilde X)\to\Mm_i^r(X)
\]
is a trivial fibration with fiber 
$S_0\Mm(\tilde\CP^2)$.
\end{teor}

The proof will be done using the analytic definition of the moduli space.
Fix $SU(r)$ bundles $E\to\tilde X$, $E_X\to X$,
$E_P\to\tilde\CP^2$ with $c_2(E)=k$, $c_2(E_X)=i$ and
$c_2(E_P)=k-i$. We will use the notation $\Mm(\tilde X,E)$,
$\Mm(X,E_X)$, $\Mm(\tilde\CP^2,E_P)$ for the moduli spaces.

We begin by introducing the enlarged moduli spaces:

\begin{defi}
let $U\subset X$ be an open topological
ball around $x_0$ intersecting $C$ 
in a non-empty disk and let
$\tilde U=\pi^{-1}(U)$. 
Then we define
\begin{enumerate}
\item $\Mm^U(X,E_X)$ is the quotient by $Aut(E_X)$ of the
space of triples $(\bar\partial_X,\phi_X,\psi_X)$ where
$\bar\partial_X$ 
is a holomorphic structure on $E_X$, $\phi_X:E_X|_{C}\to\Ee_0$
is an isomorphism and $\psi_X$ is a holomorphic trivialization of
$E|_U$ that agrees with $\phi_X$ in $U\cap C$;
\item $S_i\Mm^{\tilde U\setminus L}(\tilde X,E)$ 
is the quotient by $Aut(E)$ of the
space of triples $(\bar\partial,\phi,\psi)$ where
$\bar\partial$ 
is a holomorphic structure on $E$ such that $c_2(\pi_*E^{\vee\vee})=i$,
$\phi:E|_{C}\to\Ee_0$ is an isomorphism
and $\psi$ is a holomorphic trivialization of
$E|_{\tilde U\setminus L}$ (this bundle is always trivial)
that agrees with $\phi$ in $\tilde U\cap C$;
\item $S_0\Mm^{\tilde\CP^2\setminus L}(\tilde\CP^2,E_P)$
is the quotient by $Aut(E_P)$ of the
space of triples $(\bar\partial_P,\phi_P,\psi_P)$ where
$\bar\partial_P$ 
is a holomorphic structure on $E_P$ such that $c_2(\pi_*E_P^{\vee\vee})=0$,
$\phi_P:E|_{C}\to\Ee_0$ is an isomorphism
and $\psi_P$ is a holomorphic trivialization of
$E_P|_{\tilde\CP^2\setminus L}$ 
(trivial since $c_2(\pi_*E_P^{\vee\vee})=0$)
that agrees with $\phi_P$ in $C_\infty$.
\end{enumerate}
\end{defi}

Before we proceed we introduce the useful result:

\begin{lema}\label{lema1}
Let $B$ be a 4 dimensional ball and consider two holomorphic vector bundles
$E_1,E_2\to B$. Let $\phi:E_1|_{B\setminus0}\to E_2|_{B\setminus0}$
be an isomorphism. Then $\phi$ extends to an isomorphism
$\phi:E_1\to E_2$.
\end{lema}
\begin{proof}
$\phi$ is equivalent to a map $\phi:B\setminus0\to Gl(r,\Cc)\subset
\Cc^{r^2}$. By Hartog's theorem this map extends to a map
$\phi:B\to\Cc^{r^2}$. Composing with the determinant we get a map
$det\circ\phi:B\to\Cc$ which can only vanish at $0\in B$, hence it never
vanishes. We conclude that the image of $\phi$ lies in $Gl(r,\Cc)$.
\end{proof}

\begin{prop}
The spaces $\Mm^U(X,E_X)\times
S_0\Mm^{\tilde\CP^2\setminus L}(\tilde\CP^2,E_P)$ and \linebreak
$S_i\Mm^{\tilde U\setminus L}(\tilde X,E)$ are homeomorphic.
\end{prop}
\begin{proof}
We divide the proof into four steps:
\begin{enumerate}
\item  
The first goal is to define a map
\[
g:\Mm^U(X,E_X)\times
S_0\Mm^{\tilde\CP^2\setminus L}(\tilde\CP^2,E_P)\to
\Mm^{\tilde U\setminus L}(\tilde X,E).
\]
Let $[\bar\partial_X,\phi_X,\psi_X]\in\Mm^U(X,E_X)$,
$[\bar\partial_P,\phi_P,\psi_P]\in
S_0\Mm^{\tilde\CP^2\setminus L}(\tilde\CP^2,E_P)$.
Let
\[
E_{XP}=E_X|_{X\setminus\{x_0\}}\bigcup_{\psi_P^{-1}\psi_X}
E_P|_{\tilde U}
\]
We claim $E_{XP}$ is isomorphic to $E$ as topological vector bundles.
It is enough to show that $c_2(E_{XP})=k$. To prove this fix trivializations
$h_X$ of $E_X$ on $X\setminus\{x_0\}$ and $h_P$ of $E_P$ on $\tilde U$.
Then,
\[
E_X\simeq \Cc_{X\setminus\{x_0\}}\bigcup_{\psi_Xh_X^{-1}}\Cc_U\ ,\ 
E_P\simeq \Cc_{\tilde\CP^2\setminus L}
\bigcup_{h_P\psi_P^{-1}}\Cc_{\tilde U}
\]
Hence, seen as maps $S^3\to Gl(r,\Cc)$, $\psi_Xh_X^{-1}$ and
$h_P\psi_P^{-1}$ have degree $i$ and $k-i$ respectivelly.
Now
\[
E_{XP}\simeq \Cc_{X\setminus\{x_0\}}
\bigcup_{\psi_Xh_X^{-1}h_P\psi_P^{-1}}\Cc_{\tilde U}
\]
hence $c_2(E_{XP})=k$.

Now we define the map $g$.
$\bar\partial_X$ and $\bar\partial_P$ 
induce a holomorphic structure $\bar\partial_{XP}$ on
$E_{XP}$ and we have 
$[\bar\partial_{XP},\phi_X,\psi_X|_{\tilde U\setminus L}]\in
\Mm^{\tilde U\setminus L}(\tilde X,E_{XP})$.

Choose an isomorphism
$f:E_{XP}\to E$. Then $f$ induces a map
$f_{\#}:\Mm^{\tilde U\setminus L}(\tilde X,E_{XP})
\to\Mm^{\tilde U\setminus L}(\tilde X,E)$. We define
\[
g([\bar\partial_X,\phi_X,\psi_X],[\bar\partial_P,\phi_P,\psi_P])=
f_{\#}([\bar\partial_{XP},\phi_X,\psi_X|_{\tilde U\setminus L}])
\]
This map does not depend on the choice of isomorphism $f$.
\item
Now we show that the image of $g$ lies in
$S_i\Mm^{\tilde U\setminus L}(\tilde X,E)$.
Let $\Ee\to\tilde X$ be the bundle $E_{XP}$ with holomorphic structure
induced by $\bar\partial_X$ and $\bar\partial_P$.
We claim that $(\pi_*\Ee)^{\vee\vee}$ is biholomorphic to
$E_X$ with holomorphic structure $\bar\partial_X$.
This is a consequence of lemma \ref{lema1} since
the restriction of the bundles to $X\setminus \{x_0\}$
are biholomorphic. Hence the image of $g$ lies in
$S_i\Mm^{\tilde U\setminus L}(\tilde X,E)$.
\item Now we show $g$ is continuous. 
Fix $q_X=(\bar\partial_X,\phi_X,\psi_X)\in\Cc^U(X,E_X)$,
$q_P=(\bar\partial_P,\phi_P,\psi_P)\in
\Cc^{\tilde\CP^2\setminus L}(\tilde\CP^2,E_P)$.
Fix balls $B_r(x_0)\subset B_R(x_0)\subset U$ and let
$K=\overline{B_R(x_0)}\setminus B_r(x_0)$.
Choose $\varepsilon>0$ such that $W=B_\varepsilon(\Id)\subset Gl(r,\Cc)$ 
Then $(K,W)$ defines an open neighborhood of $\Id$ in
$C^\infty(U,Gl(r,\Cc))$. 

For each pair 
$q=([\bar\partial_X,\phi_X,\psi_X],[\bar\partial_P,\phi_P,\psi_P]))$ let
\[
E_q=E_X|_{X\setminus\{x_0\}}\bigcup_{\psi_P^{-1}\psi_X}E_P|_{\tilde U}
\]
Continuity of $g$ will follow if we construct a continuous family
of isomorphisms $f_q:E_q\to E$, or, what amounts to the same,
a continuous family $f_q:E_q\to E_{q_0}$ for fixed
\[
q_0=([\bar\partial_{X0},\phi_{X0},\psi_{X0}],
[\bar\partial_{P0},\phi_{P0},\psi_{P0}]))
\]
A map $f_q:E_q\to E_{q_0}$ is equivalent to the diagram
\[
\xymatrix{
E_X|_{X\setminus\{x_0\}}\ar[rr]^{f_X}\ar[d]^{\psi_X}&&
E_X|_{X\setminus\{x_0\}}\ar[d]^{\psi_{X0}}\\
U\setminus\{x_0\}\times\Cc^r\ar[rr]^{f_U}&&
U\setminus\{x_0\}\times\Cc^r\\
E_P|_{\tilde U}\ar[rr]^{f_P}\ar[u]^{\psi_P}&&
E_P|_{\tilde U}\ar[u]^{\psi_{P0}}
}\]
Assume $\psi_{X0}\psi_X^{-1}$ and $\psi_{P0}\psi_P^{-1}$ are
in a $(K,W)$ neighborhood of $q_0$.
Fix a monotonous $C^\infty$ function $\eta$ such that $\eta(\rho)=0$
for $\rho<r$ and $\eta(\rho)=1$ for $\rho>R$. $\eta$ induces a map
$\tilde\eta:K\to\Rr$, $x\mapsto\eta(|x|)$.
We define
\begin{align*}
& f_X=\Id && |x|>R\\
& f_U=\tilde\eta\psi_{X0}\psi_X^{-1}+(1-\eta)\psi_{P0}\psi_P^{-1} && r<|x|<R\\
& f_P=\Id && |x|<r
\end{align*}
This completes the proof of continuity of $g$.
\item
Now we construct maps
\begin{gather*}
g_X:S_i\Mm^{\tilde U\setminus L}(\tilde X,E)\to\Mm^U(X,E_X)\\
g_P:S_i\Mm^{\tilde U\setminus L}(\tilde X,E)\to
S_0\Mm^{\tilde\CP^2\setminus L}(\tilde\CP^2,E_P)
\end{gather*}
Let $[\bar\partial_{\tilde X},\phi,\psi]\in 
S_i\Mm^{\tilde U\setminus L}(\tilde X,E)$.
Define the bundles
\begin{gather*}
\hat E_X=E|_{\tilde X\setminus L}\bigcup_{\psi}U\times\Cc^r\\
\hat E_P=E|_{\tilde U}\bigcup_{\psi}\tilde\CP^2\setminus L\times\Cc^r
\end{gather*}
Then $\bar\partial_{\tilde X}$ 
induces holomorphic structures $\bar\partial_X$ on $\hat E_X$ and
$\bar\partial_P$ on $\hat E_P$. Proceeding as above we choose isomorphisms
$f_X:\hat E_X\to E_X$ and $f_P:\hat E_P\to E_P$. Then we define
\begin{gather*}
g_X([\Ee,\phi,\psi])=(f_X)_\#([\hat\Ee_X,\phi_X,\Id])\\
g_P([\Ee,\phi,\psi])=(f_P)_\#([\hat\Ee_P,\Id,\Id])
\end{gather*}
The proof that $g_X,g_P$ are well defined continuous maps
proceeds as the corresponding proof for $g$.
To conclude the proof we observe that
$g\circ(g_X\times g_P)=\Id$ and $(g_X\times g_P)\circ g=\Id$.
\end{enumerate}
\end{proof}

Now we look more closely at the map 
$g_X:S_i\Mm^{\tilde U\setminus L}(\tilde X,E)\to\Mm^U(X,E_X)$ introduced 
in the proof of the previous proposition.

\begin{prop}
Consider the projection maps
\begin{align*}
&pr_{\tilde X}:S_i\Mm^{\tilde U\setminus L}(\tilde X,E)\to
S_i\Mm(\tilde X,E)\\
&pr_X:\Mm^U(X,E_X)\to\Mm(X,E_X)
\end{align*}
$g_X$ preserves the orbits of $pr_{\tilde X}$ and $pr_X$, hence
it induces a map $\hat g_X:S_i\Mm(\tilde X,E)\to\Mm(X,E_X)$.
Furthermore, $\hat g_X=\pi_*^{\vee\vee}=\pi_\bullet$
and is determined by the following property:
\begin{itemize}
\item Let $[\Ee,\phi]\in\Mm(\tilde X,\Ee)$. Then there is a unique
holomorphic bundle $\Ee_X\to X$ such that
$\Ee|_{\tilde X\setminus L}=\Ee_X|_{X\setminus\{x_0\}}$ and we have
$\hat g_X([\Ee,\phi])=[\Ee_X,\phi]$.
\end{itemize}
\end{prop}
\begin{proof}
We begin by proving the last statement. Let $[\Ee,\phi]\in\Mm(\tilde X,\Ee)$
and define $\Ee_X=\pi_*\Ee^{\vee\vee}$. Then clearly
$\Ee|_{\tilde X\setminus L}=\Ee_X|_{X\setminus\{x_0\}}$. Uniqueness
of $\Ee_X$ then follows from lemma \ref{lema1}.

Now, from this uniqueness property and by definition of $g_X$ it follows
that $pr_X\circ g_X([\Ee,\phi,\psi])=[\Ee_X,\phi]$. This shows
$g_X$ preserves the fibers of $pr_X$ hence $\hat g_X$ is well defined.

Now it follows that $\hat g_X=\pi_\bullet$,
\end{proof}

\begin{prop}\label{prop46}
Consider the projection maps
\begin{align*}
&pr_{\tilde X}:\Mm^{\tilde U\setminus L}(\tilde X,E)\to\Mm(\tilde X,E)\\
&pr_X:\Mm_k^U(X)\to\Mm_k(X)\\
&pr_P:\Mm_k^{\tilde\CP^2\setminus L}(\tilde\CP^2)\to
\Mm_k(\tilde\CP^2)
\end{align*}
These maps are principal bundle maps with contractible fiber. 
\end{prop}
\begin{proof}
We divide the proof into three steps:
\begin{enumerate}
\item
For an open set $A$ we define
\[
G(A)=
\left\{\eta\in{\rm Hol}(A,Gl(r,\Cc))\,:\,\eta|_{A\cap C}=\Id\right\}
\]
Then, for $A=\tilde U\setminus L,\,U,\,\tilde\CP^2\setminus L$,
$G(A)$ acts on $\Mm^A(Y)$ (where $Y=\tilde X,X,\tilde\CP^2$
respectivelly) by
\[
\eta[\Ee,\phi,\psi]=[\Ee,\phi,\eta\psi]
\]
This action is free and its orbits are the fibers of $pr,pr_X,pr_P$
respectivelly.
\item
$G(U)$ is clearly contractible. From lemma \ref{lema1} it follows
that $G(\tilde U\setminus L)=G(U)$
and $G(\tilde\CP^2\setminus L)=G(\CP^2)=\{\Id\}$. 
\item
Finally we need to show the existence of local sections.
We will first build a section $s_X$ of $pr_X$.
As in the proof of proposition \ref{prop33}, for $T_X\subset \Mm(X,E_X)$
we choose a trivialization
$\psi_X:\Ff|_{T_X\times U}\to\Oo_{T_X\times U}^r$. We can choose
$\psi_X$ compatible with $\alpha_X$ on $C\cap U$.
Now, for $t=[\Ee,\phi]\in T_X$, $\psi_X$ gives a map
$\psi_{Xt}:\Ee|_U=\imath_t^*\Ff_X|_U\to \Oo_U^r$. Then
$s_X(\Ee,\phi)=(\Ee,\phi,\psi_{Xt})$.

Sections of $pr_{\tilde X}$ $p_P$ are built in a similar way:
for $T\subset\Mm(\tilde X,E)$ and $T_P\subset\Mm(\tilde\CP^2,E_P)$
we build trivializations $\psi_{\tilde X}$ and $\psi_P$ following
the same procedure as in the proof of proposition \ref{prop33}.
\end{enumerate}
\end{proof}

\begin{coro}
The spaces
$\Mm_i(X,E_X)\times 
S_0\Mm_{k-i}(\tilde\CP^2,E_P)$  and $S_i\Mm_k(\tilde X,E)$
are homotopically equivalent.
\end{coro}

We are ready to prove theorem \ref{teor51}

\begin{proof}
Fix a global section $s:\Mm(X,E_X)\to\Mm^U(X,E_X)$ of
$pr_X$. Define 
$f_s:\Mm(X,E_X)\times S_0\Mm(\tilde\CP^2,E_P)\rightarrow
S_i\Mm(\tilde X,E)$ by $f_s=pr\circ g\circ(s\times\Id)$.
It is easy to check that $f_s$ is a bijection. We need to show that its 
inverse is continuous. To that end we construct a local inverse as follows:
let $S\subset\Mm(\tilde X,E)$ be affine. Then $s$ induces a map
$s_{\tilde X}:S\to\Mm^{\tilde U\setminus L}(\tilde X,E)$ given by
the restriction of $s\circ \hat g_X$ to $\tilde U\setminus L$
(see the construction of a trivialization $\tilde\psi$
in the proof of proposition \ref{prop33}). We define
$F_s=(pr\times\Id)\circ(g_X\times g_P)\circ s_{\tilde X}$. Then
$F_s$ is the inverse of $f_s$ which concludes the proof.
\end{proof}

\subsection{Rank Stabilization}

Direct sum with a trivial bundle induces maps 
$j:\Mm_k^r(X,E)\to\Mm_k^{r'}(X,E\oplus\Oo_X^{r'-r})$
for $r'>r$. We define the rank stable moduli space by
\[
\Mm_k^\infty(X)=\lim_{\substack{\longrightarrow\\r}}\Mm_k^r
\]
The definitions of $S_i,F_i$ carry through to the rank stable situation.
We want to extend the previous results to these moduli spaces:

\begin{teor}
The spaces $S_i\Mm_k^\infty(\tilde X)$ and
$\Mm_i^\infty(X)\times S_0\Mm_{k-i}^\infty(\tilde\CP^2)$ are homotopically
equivalent.
\end{teor}
\begin{proof}
We will show that the maps $f_s,F_s$ introduced in the last section
can be defined on the direct limits and they are homotopic inverse
of each other.

We first observe that the rank stable
enlarged moduli spaces $\Mm_k^{\infty,A}(X)$
($A=U,\tilde U\setminus L,\tilde CP^2\setminus L$)
can be defined in the same way
Consider then
$f_s=pr\circ g\circ(s\times\Id)$. $pr,g$ commute with the inclusion $j$.
So we only need to show that the diagram
\[
\xymatrix{
\Mm_k^r(X)\ar[r]^-{s_r}\ar[d]^j&\Mm_k^{r,U}(X)\ar[d]^j\\
\Mm_k^{r'}(X)\ar[r]^-{s_{r'}}&\Mm_k^{r',U}(X)
}\]
is homotopy commutative. First observe that
$pr\circ j\circ s_r = pr\circ s_{r'}\circ j=j$. Hence, 
from proposition \ref{prop46} we can find, for each 
$t=[\Ee,\phi]\in\Mm_k^r(X)$, a map
$h_t\in{\rm Hol}\,(U, Gl(r',\Cc))$ such that 
$h(j \circ s_1(t))=s_2\circ j(t)$. This defines a map
$H:\Mm_k^r(X)\to {\rm Hol}\,(U, Gl(r',\Cc))$. 
Since the space of such maps is connected, $j\circ s_r$
is homotopic to $s_{r'}\circ j$. We conclude that the map $f_s$
can be defined in the direct limit. In the same way we define
the inverse map $F_s$ in the direct limit. This concludes the proof.
\end{proof}

\section{Caracterization of $S_0\overline{\Mm_k^r(\tilde\CP^2)}$ using Monads}

The objective of this section is to give a characterization of points
in $S_0\Mm_k^r(\tilde\CP^2)$ in terms of a monad description.
This result also appears in \cite{San01} and \cite{Buc02}.
We also give an explicit description of the map
$S_0\Mm(\tilde\CP^2)\to S_0\Mm^{\tilde\CP^2\setminus L}(\tilde\CP^2)$
from section 4.

We begin by sketching the monad description of the spaces
$\Mm_k^r(\CP^2)$ and $\Mm_k^r(\tilde\CP^2)$. We follow
\cite{Kin89}. See also \cite{BrSa00}.

Let $L_\infty\subset\CP^2$ be a rational curve and let $L$ be the
exceptional divisor.
Let $W,W_0,W_1$ be $k$ dimensional complex vector spaces.
Choose sections $x_1,x_2,x_3$ spanning $H^0(\mathcal O(L_\infty))$
and $y_1,y_2$ spanning
\mbox{$H^0(\mathcal O(L_\infty-L))$} so that
$x_3$ vanishes on $L_\infty$ and $x_1y_1+x_2y_2$ spans the kernel
of
\[
H^0(\mathcal O(L_\infty))\otimes H^0(\mathcal O(L_\infty-L))
\longrightarrow H^0(\mathcal O(2L_\infty-L))
\]

\subsection{The moduli space over $\CP^2$}

Let $\mathcal R$ be the space of 4-tuples $m=(a_1,a_2,b,c)$ with
$a_i\in{\rm End}(W)$, $b\in{\rm Hom}(\Cc^r,W)$, $c\in{\rm Hom}(W,\Cc^r)$,
obeying the integrability condition $[a_1,a_2]+bc=0$.
For each $m=(a_1,a_2,b,c)\in\mathcal R$ we define maps $A_m,B_m$
\[
\xymatrix{
W(-L_\infty)\ar[r]^-{A_m}&
W^{\oplus2}\oplus\Cc^n\ar[r]^-{B_m}&
W(L_\infty)
}\]
by
\[
A_m=\left[\begin{array}{c}
x_1-a_1x_3\\
x_2-a_2x_3\\
cx_3
\end{array}\right]\,,\ 
B_m=\left[\begin{array}{ccc}
-x_2+a_2x_3&x_1-a_1x_3&bx_3
\end{array}\right]
\]
Then $B_mA_m=0$. The assignement 
$m\mapsto\Ee_m={\rm Ker\,}B_m/{\rm Im\,}A_m$
induces a map $f:\mathcal R\to\overline{\Mm_k^r(\CP^2)}$.

A point $m\in\mathcal R$ is called non-degenerate if 
$A_m$ and $B_m$ have maximal
rank at every point in $\CP^2$. Then $\Ee_m$ is locally free.

\begin{teor}
Let $\mathcal M_k^r(\CP^2)$ denote the quotient of the space of
non degenerate points in $\mathcal R$ by the action of $Gl(W)$:
\[
g\cdot(a_1,a_2,b,c)=(g^{-1}a_1g,g^{-1}a_2g,
g^{-1}b,cg)
\]
Then the map $f$ induces an isomorphism 
$\mathcal M_k^r(\CP^2)\to\Mm_k^r(\CP^2)$.
\end{teor}

For a proof see \cite{Don84}, proposition 1.

\begin{teor}\label{teormonadscomplete}
Let $\overline{\mathcal M_k^r(\CP^2)}$ be the algebraic quotient 
$\mathcal R/Gl(W)$. This space is isomorphic to the 
Donaldson-Uhlenbeck completion of the moduli space of instantons over $S^4$.
\end{teor}

For a proof see \cite{DoKr90}, sections 3.3, 3.4, 3.4.4.

For future reference we sketch here how the map
from $\mathcal R/Gl(W)$ to the
Donaldson-Uhlenbeck completion of the moduli space of instantons
is constructed (see \cite{Kin89} for details):

Let $m=(a_1,a_2,b,c)\in\mathcal R$.
A subspace $V\subset W$ is called $b$-special with respect to $m$ if 
\begin{equation}
a_i(V)\subset V\ (i=1,2)\ {\rm and}\ {\rm Im}\,b\subset V\label{1}
\end{equation}
A subspace $V\subset W$ is called $c$-special with respect to $m$ if
\begin{equation}
a_i(V)\subset V\ (i=1,2)\ {\rm and}\ V\subset{\rm Ker}\,c\label{2}
\end{equation}
$m$ is called completely reducible if for every $V\subset W$
which is $b$-special or $c$-special, there is a complement
$V'\subset W$ which is $c$-special or $b$-special respectively.

\begin{prop}\label{PropMonadsSpecial}
Let $m=(a_1,a_2,b,c)\in\mathcal R$.
\begin{enumerate}
\item $m$ is non degenerate if and only if the only $b$-special subspace
is $W$ and the only $c$-special subspace is $0$;
\item For every $m$, the orbit of $m$ under $Gl(W)$ contains in its closure
a canonical completely reducible orbit and completely reducible orbits
have disjoint closures;
\item If $m$ is completely reducible then, after acting with some
$g\in Gl(W)$ we can write
\[
a_i=\begin{bmatrix}a_i^{red}&0\\0&a_i^{\Delta}\end{bmatrix}\ ,\ 
b=\begin{bmatrix}b^{red}\\0\end{bmatrix}\ ,\ 
c=\begin{bmatrix}c^{red}&0\end{bmatrix}
\]
where $(a_1^{red},a_2^{red},b^{red},c^{red})$ is non-degenerate
and the matrices $a_1^\Delta,a_2^\Delta$ can be
simultaneously diagonalized. Such a configuration is equivalent
to the following data:
\begin{itemize}
\item An irreducible integrable configuration 
$(a_1^{red},a_2^{red},b^{red},c^{red})$
corresponding to a bundle with $c_2=l\leq k$;
\item $k-l$ points in $\Cc^2=\CP^2\setminus L_\infty$ 
given by the eigenvalue pairs of
$ a_1^\Delta,a_2^\Delta$
\end{itemize}
This is precisely the Donaldson-Uhlenbeck completion.
\end{enumerate}
\end{prop}

\subsection{The moduli space over $\tilde\CP^2$}

Let $\tilde{\mathcal R}$ be the space of 5-tuples 
$\tilde m=(a_1,a_2,d,b,c)$ with
$a_i\in{\rm Hom}(W_1,W_0)$, $d\in{\rm Hom}(W_0,W_1)$,
$b\in{\rm Hom}(\Cc^r,W_0)$, $c\in{\rm Hom}(W_1,\Cc^r)$,
such that $a_1(W_1)+a_2(W_1)+b(\Cc^r)=W_0$,
obeying the integrability condition $a_1da_2-a_2da_1+bc=0$.
For each $\tilde m=(a_1,a_2,d,b,c)\in\tilde{\mathcal R}$ 
we define maps $A_{\tilde m},B_{\tilde m}$
\begin{multline*}
W_1(-L_\infty)\oplus W_0(L-L_\infty)\overset{A_{\tilde m}}{\longrightarrow}
\left(W_0\oplus W_1\right)^{\oplus2}\oplus\Cc^n
\overset{B_{\tilde m}}{\longrightarrow}\\
\rightarrow W_0(L_\infty)\oplus W_1(L_\infty-L)
\end{multline*}
by
	\[
A_{\tilde m}=\left[\begin{array}{cc}
a_1x_3&-y_2\\
x_1-da_1x_3&0\\
a_2x_3&y_1\\
x_2-da_2x_3&0\\
cx_3&0
\end{array}\right]\,,\ 
B_{\tilde m}=\left[\begin{array}{ccccc}
x_2&a_2x_3&-x_1&-a_1x_3&bx_3\\
dy_1&y_1&dy_2&y_2&0
\end{array}\right]
\]
Then $B_{\tilde m}A_{\tilde m}=0$. 
The assignement 
$\tilde m\mapsto\Ee_{\tilde m}={\rm Ker\,}B_{\tilde m}/{\rm Im\,}A_{\tilde m}$
induces a map $\tilde f:\tilde{\mathcal R}\to\overline{\Mm_k^r(\tilde\CP^2)}$.

A point $\tilde m\in\tilde{\mathcal R}$ is called non-degenerate 
if $A_{\tilde m}$ and $B_{\tilde m}$ 
have maximal rank at every point in $\tilde\CP^2$. 

\begin{teor}
Let $\mathcal M_k^r(\tilde\CP^2)$ denote the quotient of the space of
non degenerate points in $\tilde{\mathcal R}$ 
by the action of $Gl(W_0)\times Gl(W_1)$:
\[
(g_0,g_1)\cdot(a_1,a_2,b,c,d)=(g_0^{-1}a_1g_1,g_0^{-1}a_2g_1,
g_0^{-1}b,cg_1,g_1^{-1}dg_0)
\]
Then the map $\tilde f$ induces an isomorphism 
$\mathcal M_k^r(\tilde\CP^2)\to\Mm_k^r(\tilde\CP^2)$.
\end{teor}

\subsection{The theorem and its proof}

We are now in conditions to state the theorem:

\begin{teor}\label{teormonads}
Let $\tilde m=(a_1,a_2,d,b,c)\in\tilde{\mathcal R}$. Then
$\Ee_{\tilde m}\in S_0\overline{\Mm_k^r(\tilde\CP^2)}$ if and only if
$da_1$, $da_2$ 
are nilpotent and, for any sequence $i_1,\ldots,i_n\in\{1,2\}$, we have
\[
c\left(\prod_{j=1}^nda_{i_j}\right)db=0
\]
\end{teor}

We divide the proof into several propositions.

\begin{prop}
Let $\pi_\#:\tilde{\mathcal R}\to\mathcal R$ be given
by $\pi_\#(a_1,a_2,d,b,c)=(da_1,da_2,db,c)$. Let
$m=\pi_\#\tilde m$. Then 
$\Ee_{\tilde m}|_{\tilde X\setminus L}$ is isomorphic to
$\Ee_m|_{X\setminus\{x_0\}}$.
\end{prop}
\begin{proof}
Let $m=\pi_\#\tilde m$.
Fix an isomorphism $W\cong W_1$. Let $p$ be the projection
$p:W^{\oplus 2}\oplus\Cc^r\to(W_0\oplus W_1)^{\oplus 2}\oplus\Cc^r$
with kernel $W_0^{\oplus 2}$. After restricting to
$\tilde X\setminus L$ we can rescale the sections so that $y_2=-x_1$,
$y_1=x_2$.
Then a direct verification shows that, for any $\tilde m$,
$p$ induces maps
${\rm Ker}\,B_{\tilde m}\to{\rm Ker}\,B_m$ and
${\rm Im}\,A_{\tilde m}\to{\rm Im}\,A_m$. Hence we get a map
$\Ee_{\tilde m}\to \Ee_m$. It is a direct computation to check that
this map is an isomorphism.
\end{proof}

It follows from this proposition and theorem \ref{teormonadscomplete}
that, for $\tilde m\in\tilde{\mathcal R}$,
$\Ee_{\tilde m}\in S_0\overline{\Mm_k^r(\tilde\CP^2)}$ if 
and only if $\Ee_{\pi_\#\tilde m}$ is a sheaf
whose whole charge is concentrated at $[0,0,1]\in\CP^2$, that is,
$\ell\left(\Ee_{\pi_{\#\tilde m}}^{\vee\vee}/
\Ee_{\pi_{\#\tilde m}}\right)=k$.
We proceed to study this situation:

\begin{lema}
Let $m=(a_1,a_2,b,c)$ be such that
$\Ee_m$ has its 
whose whole charge concentrated at $[0,0,1]$. Then,
after a change of basis we can write
\[
a_i=\left[\begin{array}{cc}J&*\\0&J\end{array}\right]\,,\ 
b=\left[\begin{array}{c}*\\0\end{array}\right]\,,\ 
c=\left[\begin{array}{cc}0&*\end{array}\right]
\]
where $J$ represents any nilpotent matrix in the Jordan canonical form.
\end{lema}

Before we begin the proof we observe that this lemma implies one direction
of theorem \ref{teormonads}.

\begin{proof}
We begin by
proving by induction on the charge that, after acting with an element
$g\in Gl(W)$, we can write.
\[
a_i=\left[\begin{array}{cc}a_{iu}&*\\0&a_{id}\end{array}\right]\,,\ 
b=\left[\begin{array}{c}*\\0\end{array}\right]\,,\ 
c=\left[\begin{array}{cc}0&*\end{array}\right]
\]
Clearly the configuration $(a_1,a_2,b,c)$
cannot be non-degenerate hence
there is a subspace $V$ which is either $b$-special
or $c$-special (proposition \ref{PropMonadsSpecial}, 1). 
We consider both cases:
\begin{enumerate}
\item If $m$ is $b$ special, after a change of basis we can write
\[
a_i=\left[\begin{array}{cc}a_i'&*\\0&f_i\end{array}\right]\,,\ 
b=\left[\begin{array}{c}b'\\0\end{array}\right]\,,\ 
c=\left[\begin{array}{cc}c'&*\end{array}\right]
\]
So the point
$\left(\left[\begin{smallmatrix}a_i'&0\\0&f_i\end{smallmatrix}\right],
\left[\begin{smallmatrix}b'\\0\end{smallmatrix}\right],
\left[\begin{smallmatrix}c'&0\end{smallmatrix}\right]\right)$ 
is in the closure
of the orbit of $m$. It follows then from proposition 
\ref{PropMonadsSpecial} that
$m'=(a_1',a_2',b',c')$ corresponds to an ideal bundle with charge 
concentrated at $[0,0,1]$. Hence we can apply the induction hypothesis
to $m'$. We get
\[
a_i=\left[\begin{array}{ccc}
a_{iu}'&*&*\\
0&a_{id}'&*\\
0&0&f_i
\end{array}\right]\,,\ 
b=\left[\begin{array}{c}*\\0\\0\end{array}\right]\,,\ 
c=\left[\begin{array}{ccc}0&*&*\end{array}\right]
\]
which is in the desired form.
\item If $m$ is $c$-special, after a change of basis we can write
\[
a_i=\left[\begin{array}{cc}f_i&*\\0&a_i'\end{array}\right]\,,\ 
b=\left[\begin{array}{c}*\\b'\end{array}\right]\,,\ 
c=\left[\begin{array}{cc}0&c'\end{array}\right]
\]
Applying induction hypothesis to $(a_1',a_2',b',c')$ 
as in the previous case, we can write
\[
a_i=\left[\begin{array}{ccc}
f_i&*&*\\
0&a_{iu}'&*\\
0&0&a_{id}'\\
\end{array}\right]\,,\ 
b=\left[\begin{array}{c}*\\ *\\0\end{array}\right]\,,\ 
c=\left[\begin{array}{ccc}0&0&*\end{array}\right]
\]
This is in the desired form.
\end{enumerate}
Now the condition $[a_1,a_2]+bc=0$ implies
$[a_{1u},a_{2u}]=[a_{1d},a_{2d}]=0$.
So, after a change of basis we can put all these matrices in the
Jordan canonical form. Since all charge is concentrated at $[0,0,1]$,
proposition \ref{PropMonadsSpecial} implies
the eigenvalues of these matrices are all $0$.
\end{proof}

The other direction of theorem \ref{teormonads} follows directly
from the proposition

\begin{prop}
Let $m\in\mathcal R$ be such that $a_1,a_2$ are nilpotent and,
for any $n_1,n_2$ we have
\begin{equation}\label{cond.2}
ca_1^{n_1}a_2^{n_2}b=0
\end{equation}
Then $\Ee_m|_{\CP^2\setminus[0,0,1]}$ is free.
\end{prop}
\begin{proof}
We will build an explicit trivialization.
First we need to introduce some notation.
Let $\CP^2=\{[x_1,x_2,x_3]\}$ and define 
an open cover of $\CP^2\setminus[0,0,1]$ by
$U_1=\{x_1\neq0\}$, $U_2=\{x_2\neq0\}$.
Let $e_1,\ldots,e_r$ be the canonical basis of $\Cc^r$. 
Choose coordinates
$(\alpha_2,\alpha_3)\mapsto[1,\alpha_2,\alpha_3]$
in $U_1$ and
$(\beta_1,\beta_3)\mapsto[\beta_1,1,\beta_3]$
in $U_2$. Then define functions $s_i^j:U_j\to{\rm Ker\,}B_m$ by
\begin{gather*}
s_i^1=\left(0,-\alpha_3(\Id-\alpha_3 a_1)^{-1}be_i,e_i\right)\\
s_i^2=\left(\beta_3(\Id-\beta_3a_2)^{-1}be_i,0,e_i\right)
\end{gather*}
(since $a_1,a_2$ are nilpotent, $\Id-\lambda a_i$ is invertible for any
$\lambda$).
It is a direct verification that indeed the image of $s_i^1,s_i^2$
lie in $\text{Ker}\,B_m$. We want to show that $s_i^1,s_i^2$ induce
a trivialization of $\Ee_m$. We will have to show that
$S_i^1-s_i^2\in\text{Im}\,A_m$ in $U_1\cap U_2$. Then
$s^1_i,s^2_i$ induce a section of $\Ee_m$. We will show these sections are
linearly independent.
\begin{enumerate}
\item We begin by looking at the space
${\rm Im}\,A+{\rm Im}s^1_1+\ldots+{\rm Im}s^1_r$.
We can
represent the image of $s_i^1$ in terms of column vectors in matrix form.
Then, joining this matrix with $A_m$ we have
\[
\mathfrak A_1=\left[\begin{array}{cc}
\Id-\alpha_3a_1&0\\
\alpha_2-\alpha_3a_2&-\alpha_3(\Id-\alpha_3a_1)^{-1}b\\
\alpha_3c&\Id
\end{array}\right]
\]
This matrix clearly has maximum rank since $\Id-\alpha_3a_1$ is non-singular.
Hence, for dimensional reasons its columns form a basis for
${\rm Ker}\, B$. In particular $s^1_i$ are linearly independent.
\item Now we repeat the argument for $s_i^2$.
In $U_2$ we have a similar matrix:
\[
\mathfrak A_2=\left[\begin{array}{cc}
\beta_1-\beta_3a_1&\beta_3(\Id-\beta_3a_2)^{-1}b\\
\Id-\beta_3a_2&0\\
\beta_3c&\Id
\end{array}\right]
\]
which has also clearly maximum rank. Its columns form a basis for
${\rm Ker}\, B$. In particular $s^2_i$ are linearly independent.
\item
Now we will show $s^1_i-s^2_i\in{\rm Im}\,A$ in $U_1\cap U_2$.
We have, with $\beta_1=\alpha_2^{-1}$ and $\beta_3=\alpha_2^{-1}\alpha_3$,
\[
s^2(\alpha_2^{-1},\alpha_2^{-1}\alpha_3)=\left[\begin{array}{c}
\alpha_3(\alpha_2-\alpha_3a_2)^{-1}b\\0\\\Id
\end{array}\right]
\]
Since $s^2$ is in the kernel of $B$, it must be in the image of
the surjective matrix
$\mathfrak A_1$. So we can solve
\[
\left[\begin{array}{ccc}
\Id-\alpha_3a_1&0\\
\alpha_2-\alpha_3a_2&-\alpha_3(\Id-\alpha_3a_1)^{-1}b\\
\alpha_3c&\Id
\end{array}\right]
\left[\begin{array}{c}\xi_1\\\xi_2\end{array}\right]=
\left[\begin{array}{c}
\alpha_3(\alpha_2-\alpha_3a_2)^{-1}b\\0\\\Id
\end{array}\right]
\]
We obtain immediately
$\xi_1=\alpha_3(\Id-\alpha_3a_1)^{-1}(\alpha_2-\alpha_3a_2)^{-1}b$.
Now equation \ref{cond.2} implies $c\xi_1=0$. From here it follows
immediately that $\xi_2=\Id$ hence $s_i^2-s_i^1=A_m\xi_1$.
\end{enumerate}
This completes the proof.
\end{proof}

Notice that this proof gives an explicit description of the map
$S_0\Mm(\tilde\CP^2)\to S_0\Mm^{\tilde\CP^2\setminus L}(\tilde\CP^2)$.

\appendix

\section{Direct Image}

Let $S$ be a scheme and let $s\in S$. 
let $\imath_s:X\to S\times X$
and $\tilde\imath_s:\tilde X\to S\times\tilde X$
be the inclusions $x\mapsto (s,x)$.

\begin{lema}\label{appnnstep1}
For any sheaf $\Ff$ over $S\times\tilde X$ we have
\[
\imath_s^*(\Id\times\pi)_*\Ff\cong
\pi_*\tilde\imath_s^*\Ff
\]
\end{lema}
\begin{proof}
Let $V\subset X$. Then
\begin{gather*}
\imath_s^*(\Id\times\pi)_*\Ff(V)=\left(
\lim_{\substack{\displaystyle\longrightarrow\\s\times V\subset W}}
\Ff\left(\,(\Id\times\pi)^{-1}W\,\right)\right)\bigotimes_{
\underset{\substack{\longrightarrow\\
\scriptscriptstyle s\times V\subset W}}{\lim}
\Oo_{S\times X}(W)}\Oo_X(V)\\
\pi_*\tilde\imath_s^*\Ff(V)=\left(
\lim_{\substack{\displaystyle\longrightarrow\\s\times\pi^{-1}(V)
\subset U}}\Ff(U)\right)\bigotimes_{
\underset{\substack{\longrightarrow\\
\scriptscriptstyle s\times\pi^{-1}(V)\subset U}}{\lim}
\Oo_{S\times\tilde X}(U)}\Oo_{\tilde X}(\pi^{-1}(V))
\end{gather*}
which we can rewrite as
\begin{gather*}
\imath_s^*(\Id\times\pi)_*\Ff(V)=
\lim_{\substack{\displaystyle\longrightarrow\\U_1\in\mathcal S_1}}
\Ff(U_1)\bigotimes_{
\underset{\substack{\longrightarrow\\
\scriptscriptstyle U_1\in\mathcal S_1}}{\lim}
\Oo_{S\times\tilde X}(U_1)}\Oo_X(V)\\
\pi_*\tilde\imath_s^*\Ff(V)=
\lim_{\substack{\displaystyle\longrightarrow\\U_2\in\mathcal S_2}}
\Ff(U_2)\bigotimes_{
\underset{\substack{\longrightarrow\\
\scriptscriptstyle U_2\in\mathcal S_2}}{\lim}
\Oo_{S\times\tilde X}(U_2)}\Oo_X(V)
\end{gather*}
where
\begin{align*}
&\mathcal S_1=\left\{\,(\Id\times\pi)^{-1}W\,|\,s\times V\subset W\,\right\}\\
&\mathcal S_2=\left\{\,U\,|\,(\Id\times\pi)^{-1}(s\times V)\subset U\,\right\}
\end{align*}
We claim that $\mathcal S_1=\mathcal S_2$: Just observe that
if $U_2\in\mathcal S_2$ then $U_2=(1\times\pi)^{-1}(1\times\pi)(U_2)$.
This concludes the proof.
\end{proof}

\bibliography{bi}

\end{document}